%% file: sobczyk29-12-2020.tex
\documentclass[10]{article}

\input{macros}

\usepackage{array}
\usepackage{tabularx}
\usepackage{amsfonts}
\usepackage{amssymb}
\usepackage{epstopdf}
\usepackage{tikz-cd}
\usetikzlibrary{matrix,arrows,decorations.pathmorphing}

\def\bpm{\begin{pmatrix}}
\def\epm{\end{pmatrix}}
\makeindex

\title{Nested Coordinate Systems in Geometric Algebra}
\author{Garret Sobczyk \\
	Universidad de las Am\'ericas-Puebla \\
	Departamento de Actuaría F\'isica y Matem\'aticas \\
	72820 Puebla, Pue., M\'exico}
\begin{document}

\maketitle
\begin{abstract} 
A nested coordinate system is a reassigning of independent variables to take advantage of geometric or symmetry properties of a particular application. Polar, cylindrical and spherical coordinate systems are primary examples of such a regrouping that have proved their importance in the separation of variables method for solving partial differential equations. Geometric algebra offers powerful complimentary algebraic tools that are unavailable in other treatments.
 
 \smallskip
\no {\em AMS Subject Classification MSC-2020:} 15A63, 15A67, 42B37.

\smallskip
\no {\em Keywords:} Clifford algebra, coordinate systems, geometric algebra, separation of variables.
 
\end{abstract}

 \section*{0\quad Introduction}

Geometric algebra $\G_3$ is the natural generalization of Gibbs Heaviside vector algebra, but unlike the latter, it can be immediately generalized to higher dimensional geometric algebras $\G_{p,q}$ of a quadratic form. On the other hand, Clifford analysis, the generalization of Hamilton's quaternions, is also expressed in Clifford's geometric algebras \cite{AS}. The main purpose of this article is to formulate the  concept of a {\it nested coordinate system}, a generalization of the well-known methods of orthogonal coordinate systems to apply to any coordinate system. We restrict ourselves to the geometric algebra $\G_3$ because of its close relationship to the Gibbs-Heaviside vector calculus \cite{MT}. This restriction also draws attention to the clear advantages of geometric algebra over the later, because of its powerful associative algebraic structure. 

The idea of a nested rectangular coordinate system arises naturally when studying properties of polar coordinates in the $2$ and $3$-dimensional Euclidean vector spaces $\R^2$ and $\R^3$.  We begin by discussing the relationship between ordinary polar coordinates and the nested rectangular coordinate system ${\cal N}_{1,2}$, before going on to the higher dimensional nested coordinate system ${\cal N}_{1,2,3}$ utilized in the reformulation of cylindrical and spherical coordinates. A detailed discussion of the geometric algebra $\G_3$ is not given here, but results are often expressed in the closely related well-known Gibbs-Heaviside vector analysis for the benefit of the reader.   

 \section{Polar and nested coordinates systems} 
  
      Let $\G_2:=\G (\R ^2)$ be the geometric algebra of $2$-dimensional
   Euclidean space $\R^2$. An introductory treatment of the geometric algebras $\G_1$, $\G_2$ and $\G_3$ is given in \cite{Sob2019,SNF,web}.   Most important in studying the geometry of the Euclidean plane is the position vector
    \beq \bx :=\bx[x,\hat x]= x\hat x \label{position2vec} \eeq 
    expressed here as a product of its {\it Euclidean magnitude} $x$ and its {\it unit direction}, the unit vector $\hat x$.  In terms of rectangular coordinates $(x_1,x_2)\in \R^2 $, 
   \beq \bx =\bx[x_1,x_2]=x_1 e_1+x_2 e_2, \label{coorvecx} \eeq
   for the orthogonal unit vectors $e_1,e_2$ along the $x_1$ and $x_2$ axis, respectively. The advantage of our notation is that it immediately generalizes to $3$ and higher dimensional spaces of arbitrary signature $(p,q)$ in any of the definite geometric algebras $\G_{p,q}:=\G(\R^{p,q} ) $ of a quadratic form.
   
   The {\it vector derivative}, or {\it gradient} in the Euclidean plane is defined by
   \beq \nabla := e_1 \partial_1  + e_2\partial_2   \label{vecderiv} \eeq 
   where $\partial_1 := \frac{\partial }{\partial x_1}$ and $\partial_2 := \frac{\partial }{\partial x_2}$ are partial derivatives \cite[p.105]{MT}. Clearly,
   \[e_1= \partial_1 \bx = e_1 \cdot \nabla \bx, \ \ e_2= \partial_2 \bx = e_2 \cdot \nabla \bx. \]
   Since $\nabla$ is the usual 2-dimensional gradient, it has the well-known properties
   \[ \nabla \bx = 2, \quad {\rm and} \quad \nabla x = \hat x.\footnote{Note in geometric algebra, unlike in standard vector analysis, we need not write
   	$\nabla \cdot \bx =2$. This has many important consequences in the development of the subject.}  
    \]
      
      With the help of the {\it product rule} for differentiation,  
   \beq  2 =\nabla \bx = (\nabla x)\hat x + x (\nabla \hat x)= \hat x^2 + x  (\nabla \hat x). \label{vecderivhx}
    \eeq     
   Since in geometric algebra $\bx^2 = x^2$, it follows that $\hat x^2=1$, so that for $\bx \in \R^2$,
   \beq  \nabla \hat x = \frac{1}{x} \ \ {\rm and} \ \ e_1 \cdot \nabla x=e_1\cdot \hat x= \frac{x_1}{x}, \ \  e_2 \cdot \nabla x=e_2\cdot \hat x= \frac{x_2}{x} .   \label{vecDformulas} \eeq 
   Similarly, $\nabla \hat x = \frac{n-1}{x} $ for $\bx \in \R^n$. This is the first of many demonstrations of the power of geometric algebra over standard vector algebra.
  
    By a {\it nested} rectangular coordinate system ${\cal N}_{1,2}(x_1,x[x_1,x_2])$, we mean 
   \[ \bx=x \hat x =\bx [x_1,x] =\bx \big[x_1,x[x_1,x_2]\big]. \]
   The grouping of the variables allows us to consider $x_1$ and $x:=\sqrt{x_1^2+x_2^2}$ to be independent. The partial derivatives with respect to these independent variables is denoted by $\hat\partial_1 := \frac{\hat\partial}{\partial x_1}$ and  $\hat\partial_x := \frac{\hat\partial}{\partial x}$, the hat on the partial derivatives indicating the new choice of independent variables.

      For polar coordinates  $(x,\theta) \in \R^2$,  for $x:=\sqrt{x_1^2+x_2^2}\ge 0$, $0\le \theta < 2 \pi$, and $\bx :=\bx[x,\theta]$, 
    \beq \bx =x\hat x[\theta]=x(e_1 \frac{x_1}{x} + e_2 \frac{x_2}{x})= x(e_1  \cos\theta + e_2  \sin \theta),   \label{nestedpolar} \eeq
    where $\cos \theta :=\frac{x_1}{x}$ and $\sin \theta := \frac{x_2}{x}$.
 Using (\ref{vecDformulas}), 
\beq \nabla \hat x = \nabla \hat x[\theta]=(\nabla \theta)\frac{\partial \hat x}{\partial \theta} = \frac{1}{x} \quad \iff \quad  \nabla \theta = \frac{1}{x}\frac{\partial \hat x}{\partial \theta} , \ \nabla^2 \theta=0, \label{hattheta} \eeq
since 
\[ \nabla \hat x= (\nabla \theta) \partial_\theta( e_1\cos\theta + e_2 \sin \theta)=(\nabla \theta)( -e_1\sin\theta + e_2 \cos \theta ),\]
and
\[ \nabla^2 \theta = -\frac{\hat x}{x^2}\partial_{\theta} \hat x+\frac{1}{x}(\nabla \theta) \partial_{\theta}^2 \hat x = -2\Big(\frac{\hat x}{x^2}\cdot(\partial_{\theta} \hat x)\Big)=0. \]
The $\iff$ follows by multipling both sides of the first equation by the unit vector $  \partial_\theta \hat x$, which is allowable in geometric algebra. Note also the use of the famous geometric algebra identity $2\ba \cdot \bb =(\ba \bb +\bb \ba)$ for vectors $\ba$ and $\bb$, \cite[p.26]{Sob2019}. 

    The 2-dimensional gradient $\nabla$,
    \beq \nabla =e_1 \frac{\partial}{\partial x_1}+e_2\frac{\partial}{\partial x_2} =e_1 \partial_1 +e_2 \partial_2  \label{nablax}    \eeq 
    already defined in (\ref{vecderiv}), and the Laplacian $\nabla^2$ is given by
     \beq \nabla^2 = \frac{\partial^2}{\partial x_1^2}+\frac{\partial^2}{\partial x_2^2} = \partial_1^2 + \partial_2^2.  \label{nabla2x}   \eeq
   In polar coordinates,
  \beq \hat \nabla = (\nabla x)\hat\partial_x + (\nabla \theta)\hat \partial_\theta= \hat x\, \hat\partial_x + \frac{1}{x}(\hat\partial_\theta x) \hat\partial_{\theta}  \label{nablap} \eeq
  for the gradient where $\hat\partial_\theta:=\frac{\hat\partial}{\partial \theta}$,  and since $\hat \nabla^2 \theta=0$,
  \[ \hat \nabla^2 =\hat \nabla \big(  \hat x\, \hat\partial_x +(\hat\nabla \theta)\, \hat\partial_{\theta}\big) =\Big(\hat \nabla \hat x + \hat x \cdot \hat \nabla\Big) \hat \partial_{x}+\Big(\hat \nabla^2 \theta + (\hat \nabla \theta) \cdot \hat \nabla\Big) \hat \partial_{\theta} \] 
  \beq = \hat \partial_x^2 + \frac{1}{x}\hat \partial_x +\frac{1}{x^2} \hat \partial_\theta^2 . \label{nabla2p} \eeq
  for the Laplacian. The decomposition of the Laplacian (\ref{nabla2p}), directly implies that Laplace's differential equation is separable in polar coordinates.
   
   When expressed in nested rectangular coordinates ${\cal N }_{1,2}(x_1,x)$, the gradient $\nabla \equiv \hat \nabla$ takes the form
    \beq  \hat \nabla :=(\nabla x_1)\frac{\hat\partial}{\partial x_1}+(\nabla x)\frac{\hat\partial}{\partial x}=  e_1 \hat \partial_1+\hat x \hat\partial_x.    \label{nestpolar}  \eeq
   Dotting equations (\ref{nablax}) and (\ref{nestpolar}) on the left by $e_1$ and $\hat x$ gives the transformation rules 
    \[\partial_1 = \hat \partial_1 +\frac{x_1}{x}\hat \partial_x, \ \ \hat x \cdot \hat\nabla= \frac{x_1}{x} \hat \partial_1 + \hat \partial_x =\cos\theta \, \partial_1+\sin \theta \, \partial_2.\]
Using these formulas the nested Laplacian takes the form 
 \beq \hat \nabla^2 = \hat\partial_1^2+2\frac{x_1}{x} \hat\partial_x  \hat\partial_1 + \frac{1}{x}\hat\partial_x+\hat\partial_x^2=- \hat\partial_1^2 +2 \partial_1 \hat \partial_1+\frac{1}{x}\hat \partial_x+\hat \partial_x^2.  \label{templnestpolar} \eeq  
 The unusual feature of the nested Laplacian is that it is defined in terms of both the ordinary partial derivative $\partial_1$ and the nested partial derivative $\hat \partial_1$.
  Whereas partial derivatives generally commute, partial derivatives of different types do not. For example, it is easily verified that
 \[ \partial_1 \hat \partial_1 x_1 x^2=2x_1, \ \ {\rm whereas} \ \ \hat\partial_1\partial_1 x_1 x^2 = 4x_1.  \]

  Because the mixed partial derivatives $\hat \partial_x \hat \partial_1$ occurs in (\ref{templnestpolar}), Laplace's differential equation in the real rectangular coordinate system ${\cal N }_{1,2}(x_1,x)$ is not, in general, separable. Indeed, suppose that
   a harmonic function $F$ is separable, so that $F = X_1X$ 
 for $X_1=X_1[x_1], X=X[x]$. Using (\ref{templnestpolar}), 
 \beq \frac{\hat \nabla^2 F}{X_1 X}= \frac{ \partial_1^2X_1}{X_1} + \frac{\Big(\partial_x^2X+\frac{1}{x} \partial_x\Big) X}{X} 
 +  2\Big(\frac{x_1 \partial_1 X_1}{X_1}\Big)\Big(\frac{\partial_x X}{xX}\Big)=0.    \label{laplacianh} \eeq
 The last term on the prevents F in general from being separable. However, it is easily checked that $F=k \frac{x_1}{x^2}$ is harmonic and
 a solution of (\ref{templnestpolar}).
 When $X_1[x_1]=k x_1$, it is easily checked that $\frac{x_1 \partial_1 X_1}{X_1}=1 $.
 Letting $F=k x_1 X[x]$, and requiring $\hat \nabla^2F=0$, leads to the differential equation for $X[x]$,
 \[ 3\partial_x X+x \partial_x^2X =0, \]
 with the solution $X[x]=c_1 \frac{1}{x^2}+c_2$.
  The simplest example of a harmonic function $F=X_1 X$ is when $X_1=x_1$ and $X=\frac{1}{x^2}$. A graph of this function is shown in Figure 1.
   \begin{figure}[h]
  	 \begin{center}     
     	\includegraphics[width= 12.5cm,height=10.0cm,angle=0]{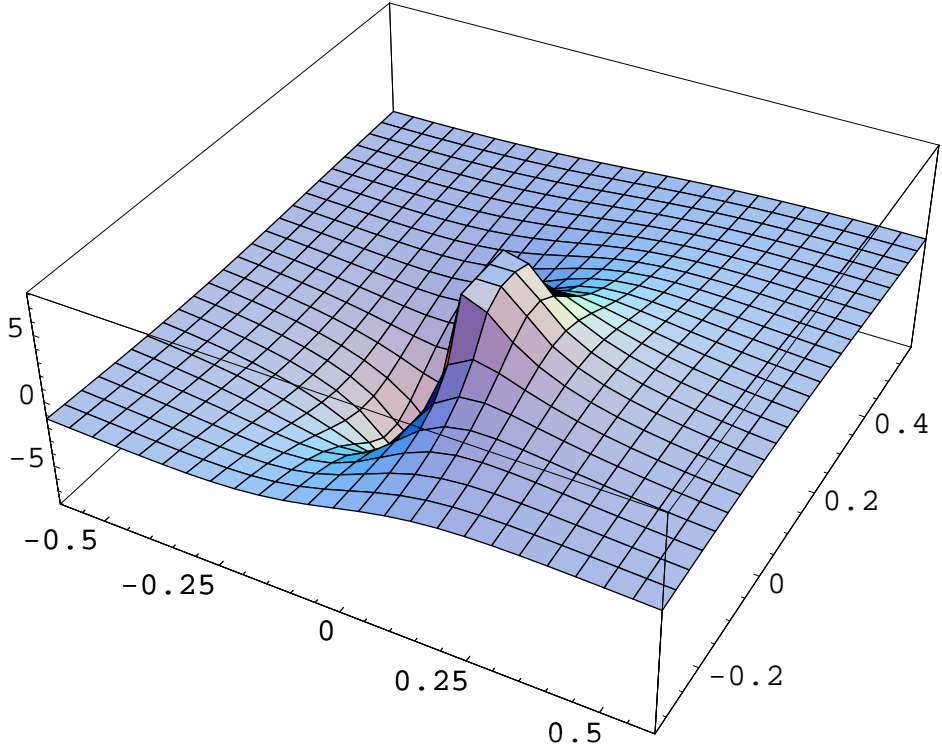}
   	\end{center}
  	\caption[harmonic]{The harmonic $2$-dimensional function $F=\frac{x_1}{x_1^2+x_2^2}$ is shown.}
  	\label{harmonicf}
\end{figure}
  
  \section{Special harmonic functions in nested coordinates}
  
  Consider the real nested rectangular coordinate system $(x_1,x_p,x)$, defined by
 \[ {\cal N}_{1,2,3}:=\{(x_1,x_p,x)| \ \bx = x \hat x = x_1 e_1 + x_p \hat x_p + x \hat x     \},  \]
 where $x_p=\sqrt{x_1^2+x_2^2}\ge 0, \ x=\sqrt{x_1^2+x_2^2+x_3^2}\ge 0$.
 In nested coordinates, the gradient $\nabla=e_1 \partial_1 + e_2\partial_2+e_3 \partial_3 $ takes the form
 \beq \hat\nabla= (\hat\nabla x_1) \hat \partial_1 +  (\hat\nabla x_p) \hat\partial_{p}+ (\hat\nabla x) \hat\partial_x = e_1 \hat \partial_1 + \hat x_p \hat\partial_{p}+\hat x \hat\partial_x,  \label{nestnabla3}   \eeq
 where $\hat \partial_p:=\frac{\hat \partial}{\partial x_p}$.
Formulas relating the gradients $\nabla$ and $\hat \nabla$ easily follow:
 \beq \partial_1 =\hat\partial_1 + e_1 \cdot \hat x_p\, \hat\partial_p+e_1\cdot\, \hat x \,\hat \partial_x= \hat\partial_1 +\frac{x_1}{x_p}\, \hat\partial_p+\frac{x_1}{x} \,\hat \partial_x    \label{lpartial1} \eeq
    \beq \partial_2 =  e_2 \cdot \hat x_p\, \hat\partial_p+e_2\cdot\, \hat x \,\hat \partial_x= \frac{x_2}{x_p}\, \hat\partial_p+\frac{x_2}{x} \,\hat \partial_x    \label{lpartial2} \eeq
 and
  \beq \partial_3 =e_3\cdot\, \hat x \,\hat \partial_x= \frac{x_3}{x} \,\hat \partial_x .   \label{lpartial3} \eeq
   
  For the Laplacian $\nabla^2 $ in nested coordinates, with the help of (\ref{nestnabla3}),
  \[ \hat\nabla^2 =\hat \nabla (e_1 \hat \partial_1 + \hat x_p \partial_{x_p}+\hat x \hat\partial_x) =e_1 \cdot \hat \nabla \,\hat \partial_1 + \hat \nabla \cdot \hat x_p\, \hat \partial_p + \hat \nabla \cdot \hat x \, \hat \partial_x  \]
  \[ =e_1 \cdot \hat \nabla \,\hat \partial_1 + \hat \nabla \cdot \hat x_p\, \hat \partial_p + \hat \nabla \cdot \hat x \, \hat \partial_x \] 
   
 \[ =\Big( \hat\partial_1 +\frac{x_1}{x_p}\, \hat\partial_p+\frac{x_1}{x} \,\hat \partial_x \Big) \hat \partial_1 +\Big(\frac{1}{x_p}+\frac{x_1}{x_p}\hat \partial_1+\hat \partial_p+\frac{x_p}{x}\hat \partial_p \Big)\hat \partial_p\] \[+\Big(\frac{2}{x}+ \frac{x_1}{x}\hat \partial_1+\frac{x_p}{x}\hat \partial_p+\hat \partial_x\Big) \hat \partial_x    \]
  \beq =\hat \partial_1^2+ \hat\partial_p^2+\hat\partial_x^2 +2\bigg( \frac{x_1}{x_p}\hat \partial_1 \hat \partial_{p}+ \frac{x_1}{x}\hat \partial_1 \hat \partial_{x}+ \frac{x_p}{x}\hat \partial_p \hat \partial_{x}\bigg)+\frac{1}{x_p}\hat \partial_p + \frac{2}{x}\hat \partial_x. \label{nestlaplace}  \eeq 
  Another expression for the Laplacian in mixed coordinates is obtained with the help of (\ref{lpartial1}),
  \beq \hat \nabla^2=- \hat \partial_1^2+ \hat\partial_p^2+\hat\partial_x^2 +2\bigg( \partial_1\hat \partial_1 + \frac{x_p}{x}\hat \partial_p \hat \partial_{x}\bigg)+\frac{1}{x_p}\hat \partial_p + \frac{2}{x}\hat \partial_x.      
             \label{lmixed} \eeq
 
  Suppose $F=F[x_1,x_p,x]$. In order for $F$ to be harmonic, $\hat\nabla^2F=0$. Assuming that $F$ is separable, $F=X_1[x_1]X_p[x_p]X_x[x]$, and applying the Laplacian (\ref{lmixed}) to $F$ gives
  \[ \hat\nabla^2F =(\hat \partial_1^2X_1)X_pX_x+X_1\Big( \big(\hat\partial_p^2+\frac{1}{x_p }\hat \partial_p\big)X_p\Big) X_x +  X_1X_p\Big(\frac{2}{x}\hat\partial_{x}X_x\Big)\] 
  \beq +2\bigg( \big( x_p\partial_p X_p\big)\big(\frac{1}{x} \hat \partial_{x}X_x\big)X_1    +\big( \partial_1 X_p\big) X_x+X_p \big(\partial_1 X_x\big) \bigg).  \label{lmixed2} \eeq
  We now calculate the interesting expression
  \[    \frac{ \big( x_p\partial_p X_p\big)\big(\frac{1}{x} \hat \partial_{x}X_x\big)X_1    +\big( \partial_1 X_p\big) X_x+X_p \big(\partial_1 X_x\big)   }{X_1X_pX_x} \] 
  \[=\Big(x_p \big(\partial_p\log X_p\big) \Big)\Big(\big(\frac{1}{x}\partial_x\log X_x\big)\Big) + \frac{\partial_1 \log(X_pX_x)}{X_1} . \]
       
        In general, because of the last term in (\ref{lmixed2}), a function
    $F=X_1 X_p X_x$ will not be separable. However, just as in the two dimensional case, there are $3$-dimensional harmonic solutions of the form
    $F= x_1^k x_p^m x^n$. Taking the Laplacian (\ref{nestlaplace}) of $F$, with the help of \cite{wofram}, gives
    \[ \hat \nabla^2F = (2km+m^2)x^nx_1^kx_p^{m-2}+(-k+k^2)x^nx_1^{k-2}x_p^m\]
    \[ +(2kn+2mn+n(1+n))x^{n-2}x_1^kx_p^m =0.    \]
    This last expression vanishes when the system of three equations,
    \[ \{ 2km+m^2=0, \  \ -k+k^2=0, \ \ {\rm and} \ \ 2kn+2mn+n(1+n)=0\}.     \]
    All of the distinct non-trivial harmonic solutions $F=x_1^k x_p^m x^n$ are listed in the following Table
   \begin{center}
   \beq	\begin{tabular}{|c|c|c|}
        \hline  	 
   	 	k & m & n \\
   		\hline \hline
   		1 & 0 & 0 \\
        \hline  	
   		0 & 0 & -1 \\
   		\hline
   		1 & -2 & 0 \\
   		\hline
   		1 & 0 & -3 \\
   		\hline
   		1 & -2 & 1 \\
        \hline	
      	\end{tabular} \eeq
   \end{center}
   
     \section{Cylindrical and spherical coordinates}    
  
\begin{figure}[h]
	\begin{center}     
		\includegraphics[width= 05.5cm,height=7.0cm,angle=0]{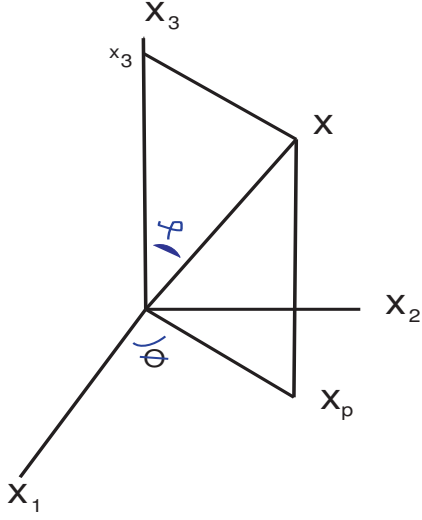}
	\end{center}
	\caption[cylsphere]{For cylindrical coordinates, $\bx = x_p \hat x_p[\theta]+x_3e_3$. For spherical coordinates, $\bx = x(e_3 \cos \varphi +\hat x_p[\theta]\sin \varphi) $. }
	\label{cylsphere}
\end{figure}

Cylindrical and spherical coordinates are examples of nested coordinates ${\cal N}_{1,2}(\R)$, and ${\cal N}_{2,3}(\R)$, respectively. For the first,
\beq \bx =\bx[x_p,\theta,x_3]= \bx_p[x_p,\theta]+\bx_3[x_3],   \label{posvecyl}\eeq
where $\bx_p=x_p\hat x_p[\theta]$, $x_p=\sqrt{x_1^2+x_2^2}$, and $\bx_3=x_3e_3$. Cylindrical coordinates $(x_p,\theta, x_3)\in \R^3=\R^2\cup \R^1$ are a decomposition of $\R^3$ into the polar coordinates $(x_p,\theta) \in\R^2$, already studied in Section 1, and $x_3 \in\R^1$. For spherical coordinates, $\bx_p =x_p  \hat x_p[\theta]$ the same as in cylindrical and polar coordinates, and 
\beq  \bx =\bx[x,\theta,\varphi]=x \hat x[\theta,\varphi])=x\Big(  e_3 \cos \varphi + \hat x_p[\theta] \sin \varphi\Big),   \label{posvecsph} \eeq
where 
\[x=\sqrt{x_1^2+x_2^2+x_3^2}, \  \hat x[\theta,\varphi]=e_3 \cos \varphi + \hat x_p[\theta] \sin \varphi, \  \hat x_p[\theta]=e_1\cos \theta + e_2 \sin \theta . \]
The basic quantities that define both cylindrical and spherical coordinates are shown in Figure \ref{cylsphere}.

The gradient $\hat \nabla$ and Laplacian $\hat \nabla^2$ for cylindrical coordinates are easily calculated. With the help of (\ref{hattheta}), (\ref{nablap}), and (\ref{nabla2p}),  
\[ \hat \nabla =(\hat \nabla x_p)  \hat\partial_p + (\hat \nabla \theta)\hat\partial_\theta + (\hat \nabla x_3)\hat\partial_3         \]
for the cylindrical gradient, and
    \beq  \hat \nabla^2 =\hat \nabla \Big(  \hat x_p\, \hat\partial_p +(\hat\nabla \theta)\, \hat\partial_{\theta}+e_3 \hat \partial_3 \Big)  = \hat \partial_p^2 + \frac{1}{x_p}\hat \partial_p +\frac{1}{x_p^2} \hat \partial_\theta^2+\hat \partial_3^2   \label{laplaceC} \eeq    
for the cylindrical Laplacian. Letting $F[\bx]=X_p[x]X_\theta[\theta]X_3[x_3]$, the resulting equation is easily separated and solved by standard methods, resulting in three second order differential equations with solutions,
\[ X_p[x_p]=k_1 J_n[\beta x_p]+k_2 Y_n[\beta x_p,]    \]
\[ X_\theta[\theta]=k_3 \cos n \theta +k_4 \sin n \theta   \]
\[ X_3[x_3]=k_5 \cosh(\alpha(m- x_3))+k_6\sinh(\alpha(m-x_3)),   \]
where $J_n$ and $Y_n$ are Bessel functions of the first and second kind . The constants are determined by the various boundary conditions that must be satisfied in different applications \cite[p.254]{TM}.

Turning to spherical coordinates $(x,\theta,\varphi)\in \R^3$, the spherical gradient
\beq \hat \nabla = (\hat \nabla x)\hat \partial_x +(\hat \nabla \theta)\hat\partial_{\theta} + (\hat \nabla \varphi)\hat\partial_\varphi  = \hat x\hat \partial_x +\frac{1}{x_p}(\hat\partial_{\theta}\hat x_p)\hat \partial_{\theta}  +\frac{1}{x}(\hat\partial_\varphi \hat x) \hat \partial_{\varphi} , \label{sphere2laplace} \eeq 
where from previous calculations for polar and cylindrical coordinates,
\beq \hat\nabla \theta = \frac{1}{x_p}(\hat \partial_p \hat x_p), \ (\hat \nabla \theta)^2 = \frac{1}{x_p^2}, \ \ \hat \nabla^2 \theta = 0, \ \ \hat \nabla \varphi=\frac{1}{x}\hat\partial_\varphi \hat x, \ \ (\hat \nabla \varphi)^2 = \frac{1}{x^2}.  \label{importprop} \eeq  
Furthermore, since $\hat x = \hat x[\theta,\varphi]=e_3\cos \varphi +\hat x_p[\theta]\sin \varphi $
\[ \frac{2}{x}=\hat \nabla \hat x =(\hat \nabla \theta)(\hat\partial_{\theta}\hat x) + (\hat \nabla \varphi)(\hat\partial_{\varphi}\hat x) =\frac{1}{x_p}(\hat\partial_{\theta}\hat x)(\hat\partial_\varphi \hat x) +\frac{1}{x} ,  \]
 it follows that
\[ (\hat\partial_{\theta}\hat x)(\hat\partial_\varphi \hat x) =\frac{x_p}{x}=\sin \varphi, \ \ {\rm and} \ \ \hat \nabla^2 \varphi =\frac{x_3}{x^2 x_p}. \]

That $ \hat \nabla^2 \varphi =\frac{x_3}{x^2 x_p}$ follows using (\ref{sphere2laplace}) and (\ref{importprop}), 
\[ \hat \nabla^2\varphi = \hat \nabla  \big(\frac{1}{x}\hat \partial_\varphi \hat x\big) = \Big(-\frac{\hat x}{x^2}+ \frac{1}{x}\hat \nabla \Big) \hat\partial_\varphi \hat x \]
\[ = -\frac{\hat x}{x^2}\hat \partial_{\varphi}\hat x-\frac{1}{x}\Big( ( \hat \nabla x) \hat \partial_{x}\hat\partial_\varphi \hat x + (\hat \nabla \theta)\hat\partial_\theta \hat \partial_\varphi\hat x+ (\hat \nabla \varphi)\hat\partial_\varphi^2 \hat x\Big)             \]
\[  =-\Big(\frac{\hat x}{x^2}\hat \partial_{\varphi}\hat x + \frac{1}{x^2}(\hat \partial_{\varphi}\hat x )\hat x\Big)+ \frac{1}{xx_p}(\hat\partial_\theta \hat x_p) (\hat \partial_\varphi \hat\partial_\theta x)=\frac{x_3}{x^2x_p},  \]
since partial derivatives commute, $\hat \partial_{x} \hat x = 0$, and $ \hat\partial_\varphi^2 \hat x = -\hat x$.

For the spherical Laplacian, using (\ref{sphere2laplace}) and (\ref{importprop}),
\[ \hat \nabla^2 =\hat \nabla \Big( \hat x\hat \partial_x +(\hat \nabla \theta)\hat \partial_{\theta}  +(\hat \nabla \varphi) \hat \partial_{\varphi}     \Big)    \]
\[ = \Big(\frac{2}{x}+   \hat \partial_x\Big)\hat \partial_x + (\hat \nabla \theta)\cdot \hat \nabla \hat \partial_{\theta}  +\Big(\hat\nabla^2 \varphi +(\hat \nabla \varphi)\cdot \hat \nabla \Big)\hat \partial_{\varphi}\]
\[ =\Big(\hat \partial_x +\frac{2}{x} \Big)\hat \partial_x  +\frac{1}{x_p^2}\hat \partial_{\theta}^2+\Big(\frac{x_3}{x^2 x_p} +\frac{1}{x^2}\hat \partial_{\varphi} \Big)  \hat \partial_{\varphi},    \]
equivalent to the usual expression for the Laplacian in spherical coordinates \cite[p.256]{TM}. 

Just as in cylindrical coordinates, the solution of Laplace's equation in spherical coordinates is separable, $F=X_x[x]X_\theta[\theta]X_\varphi[\varphi]$,
resulting in three second order differential equations with solutions
\[ X_x[x]=k_1 x^{\beta}+k_2 x^{-(\beta +1)},   \]
\[ X_\theta[\theta]=k_3 \cos n \theta +k_4 \sin n \theta ,  \]
\[ X_\varphi[\varphi]=k_5 P_n^m(\cos \varphi)+k_6 Q_n^m(\cos \varphi),   \]
where $P_n^m$ and $Q_n^m$ are the Legendre functions of the first and second kind, respectively \cite[p.258]{TM}.

\section*{Acknowledgment}
  This work was largely inspired by a current project that author has with Professor Joao Morais of Instituto Tecnol\'ogico Aut\'onomo de M\'exico, utilizing spheroidal coordinate systems. The struggle with this orthogonal coordinate system \cite{Hob1931}, led the author to re-examine the foundations of general coordinate systems in geometric algebra \cite[p.63]{SNF}.

  \end{document}

%% file: macros.tex
\def\ba{\mathbf{a}}   
\def\bb{\mathbf{b}}   
\def\bx{\mathbf{x}}   
\def\G{\mathbb{G}} 

\def\R{\mathbb{R}}  

\def\no{\noindent}

\def\beq{\begin{equation}}
\def\eeq{\end{equation}}

